\documentclass[12pt,leqno]{article}
\newcommand{\R}{{\bf R}}
\def \R{{\rm I\!\!\!R}}

\def \AR{{\cal A}}
\def \BR{{\cal B}}

\def \FR{{\cal F}}

\def \LR{{\cal L}}

\def \UR{{\cal U}}
\def \vf{\varphi}

\newcommand{\cq}{\hfill$\diamond$\\}
\newcommand{\dem}{\noindent \underline{Proof} :}
\newcommand{\ba}{\[\begin{array}{rl}}
\newcommand{\ea}{\end{array}\]}
\newcommand{\bac}{\[\left\{\begin{array}{rl}}
\newcommand{\eac}{\end{array}\right.\]}
\newcommand{\be}{\begin{equation}\label}
\newcommand{\ce}{\begin{array}{rl}}
\newcommand{\ee}{\end{array}\end{equation}}

\newcommand{\xxcu}{X^{x,\mbox{\tiny u}}}

\newcommand{\xxuv}{x^{x,u,v}}
\newcommand{\xxu}{X^{x,u}}

\def\i{{\bf 1}\hskip-2.5pt{\rm l}}

\newcommand{\argmax}{\mbox{Argmax}}

\newtheorem{theorem}{Theorem}[section]

\newtheorem{lemma}{Lemma}[section]

\newtheorem{remark}{Remark}[section]

\newtheorem{definition}{Definition}[section]

\begin{document}
\title{Another proof for the equivalence between invariance of closed sets with respect to stochastic and deterministic systems}
\bigskip
\author{Rainer Buckdahn$^*$, Marc Quincampoix$^*$, Catherine Rainer\footnote{Universit\'e de Bretagne Occidentale, Laboratoire de Math\'ematiques, Unit\'e CNRS UMR 6205, Brest, France}, Josef Teichmann\footnote{University of Technology Vienna, Institute of mathematical methods in Economics, Vienna, Austria, jteichma@fam.tuwien.ac.at. The authors gratefully acknowledge the support from the RTN network HPRN-CT-2002-00281 (European Union) and from the FWF-grant Y 328 (Austrian Science Funds)}}
\maketitle
\section*{Abstract}
We provide a short and elementary proof for the recently proved result by G.~da Prato and H.~Frankowska that -- under minimal assumptions -- a closed set is invariant with respect to a stochastic control system if and only if it is invariant with respect to the (associated) deterministic control system.

\noindent
\textit{Keywords:} invariance of closed sets, Stratonovich drift, stochastic control system, deterministic control system, stochastic Taylor expansion.

\noindent 
\textit{MSC2000:} 93E03, 60H10.
\section{Introduction}
We deal in this note with invariance of controlled stochastic differential systems. We consider a non-empty, closed subset $ K \subset \R^n $ and ask for characterizations of invariance of $ K $ with respect to a controlled stochastic differential system
\begin{equation}
\label{sds}
 \begin{array}{ll}
dX_t=b(X_t, u_t)dt+\sigma(X_t,u_t)dW_t,& t\geq 0, \\
 X_0=x\in\R^n,
\end{array}
\end{equation}
driven by a $d$-dimensional Brownian motion $W$.\\
Invariance of $ K $ here means that $P[X^{x,u}_t\in K]=1$, for all $ x \in K $, all times $ t \geq 0 $ and all admissible control processes  $ u $.\\
There exist already a lot of literature concerning invariance as well as the connected notion of viability; characterizations of both have been expressed through stochastic tangent cones 
(\cite{ada}, \cite{gt}), viscosity solutions of second-order partial differential equations (e.g. see \cite{bg},\cite{bj},\cite{bpqr},\cite{bcq},\cite{qr}) or other approaches (e.g. see \cite{dpf},\cite{m}).\\
A natural approach of the notion of invariance is to look at the associated controlled ordinary differential system:
\begin{equation}
\label{ods0}
\begin{array}{ll}
x'(t)=\tilde b(x(t),u(t))
+\sigma(x(t),u(t))v(t),& t\geq 0,\\
x(0)=x,
\end{array}
\end{equation}
where $\tilde b(x,u)$ denotes the Stratonovich drift $\tilde b(x,u)=b(x,u)-\frac 12 \sum_{i=1}^d\langle D_x\sigma^i(x,u),\sigma^i(x,u)\rangle$ and $v\in L^1_{loc}([0,\infty),\R^n)$. 
For the case without control it is well known that invariance with respect to (\ref{sds}) is equivalent to invariance with respect to the ordinary differential system (\ref{ods0}) (see \cite{d},\cite{sus}, \cite{ad}).
 Recently G.~da Prato and H.~Frankowska \cite{df} proved the result on the equivalence for controlled deterministic and stochastic systems under minimal assumptions on the involved parameters. Our aim here is to provide a {\em new, short and very elementary proof} of this intuitive equivalence result.\\
The intuition behind our main result stems from the local asymptotics of the stochastic systems, which correspond precisely to those of the deterministic system. Reading this insight, which is well-known in numerical analysis for the given stochastic differential system, in the correct way, leads us to the proof. A central step in our investigation is to show that stochastic as well as deterministic invariance is equivalent to invariance with respect to constant controls. This permits us to pass from deterministic invariance to stochastic invariance by the classical Wong-Zakai-approach to martingale problems (which can be seen as a sort of Euler-Mayurama-scheme, too). Concerning the other direction of the proof, the necessary conditions on the parameters  follow naturally from a stochastic Taylor expansion.\\
As a crucial tool we apply optimization theory, since both invariance problems can be associated with problems of minimal distance to $ K $. Hence we can also assert an equivalence between first and second order Hamilton-Jacobi-Bellman systems.

\section{Main Theorem}
\noindent Let $U$ be some compact metric space and, for $d,n\geq
1$, let $b$ be a bounded and continuous map from $\R^n\times U$ to
$\R^n$ , Lipschitz in $x\in\R^n$ uniformly in $u\in U$, and
$\sigma$ a continuous map from $\R^n\times U$ to $\R^{n\times d}$, 
differentiable with respect to $x$, such that $\sigma$ and $D_x\sigma$ are bounded and Lipschitz, both uniformly in $u$.\\ 
Let $W$ be a $d$-dimensional Brownian motion on some probability
space $(\Omega,\FR,P)$ and $(\FR_t,t\geq 0)$ the filtration
generated by $W$, satisfying the usual assumptions. We denote by
$\UR$ the set of all $U$-valued processes $(u_t)$ that are
progressively measurable w.r.t. $(\FR_t,t\geq 0)$.\\
For $(u_t)\in\UR$, we consider the controlled stochastic differential system :
\begin{equation}
\label{sde}
 \begin{array}{ll}
dX_t=b(X_t, u_t)dt+\sigma(X_t,u_t)dW_t,& t\geq 0, \\
 X_0=x\in\R^n.
\end{array}
\end{equation}
It is well known that under the above assumptions on $b$ and $\sigma$,
the system (\ref{sde}) has a unique strong solution, which we denote by $X^{x,u}$.\\
We associate to this system
the usual second order operator:
for $\varphi\in C^2(\R^n,\R)$, $x\in\R^n$ and $u\in U$,
\[ \LR_{x,u} \varphi=\langle b(x,u),\vf(x)\rangle+
\frac 12\mbox{tr}(D^2\varphi(x)\sigma(x,u)\sigma^*(x,u)).\]
We denote by $\tilde b$ the Stratonovich drift
\[\tilde b(x,u)=b(x,u)-\frac 12 \sum_{i=1}^d\langle D_x\sigma^i(x,u),\sigma^i(x,u)\rangle,\]
where $\sigma^i(x,u)$ is the $i$-th column of the matrix
$\sigma(x,u)$.\\

\noindent Furthermore we consider a non empty closed set $K\subset\R^n$.
The notion of invariance of $K$ with respect to (\ref{sde}) is defined as follows:
\begin{definition}
We say that $K$ is {\em  invariant with respect to (\ref{sde})} if, for all $x\in K$, $u\in\UR$,
and $t\geq 0$,
$P[X^{x,u}_t\in K]=1$.
\end{definition}

\noindent We also shall introduce the deterministic system :
\begin{equation}
\label{de}
\begin{array}{ll}
x'(t)=\tilde b(x(t),u(t))
+\sigma(x(t),u(t))v(t),& t\geq 0\\
x(0)=x,
\end{array}
\end{equation}
driven by the deterministic control process
$v(t)\in {\cal B}:=L^1_{loc}([0,+\infty),\R^d)$ and $u(t)\in {\cal A}:=L^{\infty}([0,\infty),U)$.
For given $x\in\R^n$, $u\in \AR$ and $v\in {\cal B}$, the solution of  (\ref{de})  will be denoted by $x^{x,u,v}$.
The associated first order operator
is, for $\varphi\in C^1(\R^n)$,
\[\LR'_{x,u}\varphi=\langle \tilde b(x,u),D\vf(x)\rangle.\]

\begin{definition}
We say that a closed set $K $ is invariant w.r.t.
(\ref{de}) if, for every $x\in K$, $u\in \AR$ and  $v\in \BR$,
$x^{x,u,v}(t)\in K$ for every $t\geq 0$.
\end{definition}

For $\varphi:\R^n\to\R$,
we denote by $\argmax_K\varphi$ the set of $x\in K$ such that $\varphi$ attains
a maximum at $x$ in $K$.\\

\noindent Our main result is the following theorem.

\begin{theorem}
The following assertions are equivalent :\\
a)
$K$ is invariant with respect to (\ref{sde});\\
b) For all $\varphi\in C^2$ and $x\in \argmax_K\varphi$, it holds
that
\begin{equation}\label{5}
\left\{
\begin{array}{l}
\sup_{u\in U}\LR_{x,u}\vf\leq 0, \\
\langle\sigma^i(x,u),D\vf(x)\rangle=0, \forall i\in\{ 1,\ldots,d\}, \forall u\in U;
\end{array}
\right.
\end{equation}
c) For all $\varphi\in C^2$ and $x\in \argmax_K\varphi$, it holds
that
\begin{equation}
\label{hjbdet}
\left\{
\begin{array}{l}
\sup_{u\in U}\LR'_{x,u}\varphi(x)\leq 0,\\
\\
\langle \sigma^j(x,u),D\vf(x)\rangle=0, \forall i\in\{ 1,\ldots,d\},\forall u\in U,\\
\\
\mbox{the matrix }A_{\vf,x}=(a_{ij})
\mbox{ with }
a_{ij}=\langle \sigma^i(x,u),D_x\langle\sigma^j(\cdot,u),D\varphi(\cdot)\rangle(x)\rangle\\
\mbox{ is symmetric and semidefinite negative;}
\end{array}
\right.
\end{equation}
d) $K$ is invariant with respect to (\ref{de});\\
e)
 For all $\varphi\in C^2$ and $x\in \argmax_K\varphi$, it holds that
\begin{equation}
\label{iiic} \sup_{u\in \,U,\; v \in \, \R ^d }\{
\LR'_{x,u}\varphi(x) + \langle\sigma(x,u)v,D\varphi(x) \rangle\} \leq 0.
\end{equation}

\end{theorem}

\begin{remark}
Applying the notations from differential geometry $\tilde b_u
\varphi(x):=\langle \tilde b(x,u),D\vf(x)\rangle$ 
and $\sigma_u^i\varphi(x):=\langle\sigma^i(x,u),D\vf(x)\rangle$, condition (\ref{hjbdet})
can be rewritten as follows:\\
For all $u\in U$, it holds that
\[\begin{array}{l} 
\tilde b_u\vf(x)\leq 0,\\
\sigma^i_u\vf(x)=0, \forall i\in\{ 1,\ldots,d\}, \\
A_{\vf,x}=(\sigma^i_u\sigma^ j_u\varphi(x))_{ij}\mbox{ 
is symmetric semidefinite negative}.\end{array}
\]
\end{remark}

\noindent The following Lemma is crucial in the proof of the Theorem.

\begin{lemma}
\label{lemmaabc}
Let  $(W_t)_{t\geq 0}$ be a standard $\;\R^d$-valued Brownian motion issued from 0 and
$(R_t)_{t\geq 0}$ a real stochastic process satisfying
 $\lim_{t\searrow 0}\frac{R_t}t=0$ in probability.\\
Let $(\alpha_i,1\leq i\leq d)\in\R^d$, $(\beta_i,1\leq i\leq d)\in\R^d$,
$( \gamma_{ij},(i,j)\in\{1,\ldots,d\}^2, i\neq j)\in\R^{d^2-d}$ and $\delta\in\R$.
Suppose that, for all $t\geq 0$, $P$-a.s.,
\begin{equation}
\label{abc} \sum_{i=1}^d
\alpha_iW^i_t+\sum_{i=1}^d\beta_i(W^{i}_t)^2+ \sum_{1\leq i\neq
j\leq d}\gamma_{ij}\int_0^tW^i_sdW^j_s+\delta t+R_t\leq 0.
\end{equation}
Then it holds that\\
i) $\alpha_i=0$, for all $i\in\{ 1,\ldots, d\} $;\\
ii) the matrix $A\in\R^{d\times d}$ defined by
\[\left\{ \begin{array}{rcl}
A_{ij}&=&\gamma_{ij},\mbox{ for } (i,j)\in\{ 1,\ldots,d\}^2,\mbox{ with } i\neq j,\\
A_{ii}&=&2\beta_i, i\in\{ 1,\ldots, d\},
\end{array}\right.\]
is symmetric and semidefinite negative;\\
iii) $\delta\leq 0$.

\end{lemma}

\dem\
 It is easy to see that
\[  \sum_{i=1}^d\beta_i\frac{(W^{i}_t)^2}{\sqrt{t}}+
\sum_{i\neq j }\gamma_{ij}\frac{\int_0^tW^i_sdW^j_s}{\sqrt{t}}+\delta\sqrt{t}+
\frac{R_t}{\sqrt{t}}\stackrel{P}{\to} 0,{\mbox as }\; t\searrow 0,\]
while,
\[\forall t\geq 0,\;\sum_{i=1}^d \alpha_i\frac{W^i_t}{\sqrt{t}}\stackrel{(d)}
{=}\sum_{i=1}^d \alpha_iW^i_1.\]
It follows that the left hand term of (\ref{abc}) divided by $\sqrt{t}$, say
$L_t$,   converges in distribution to
$\sum_{i=1}^d \alpha_iW^i_1$.
Now the assumption $P[L_t\leq 0]=1$ for all $t>0$ implies that $P[\sum_{i=1}^d \alpha_iW^i_1\leq 0]=1$, too.
It follows that, necessarily $\alpha_1=\ldots =\alpha_d=0$.\\

\noindent Using again the scaling property of Brownian motion, we have,
for all $(i,j)\in\{ 1,\ldots,d\}^2$ and for all $t\geq 0$, $\frac
1t\int_0^tW^i_sdW^j_s\stackrel{(d)}{=}\int_0^1W^i_sdW^j_s$. By
the same arguments as above, we can deduce from (\ref{abc}) that, $P$-a.s.,
\begin{equation}
\label{bgd}
 \sum_{i=1}^d\beta_i(W^{i}_1)^2+ \sum_{i\neq j}\gamma_{ij}\int_0^1W^i_sdW^j_s+\delta \leq 0.
\end{equation}
Let us focus now on a fixed  arbitrary couple of indexes $(i,j)$ with $i\neq j$. After conditioning by $\sigma(W^i_s,W^j_s,s\geq 0)$, we get from
(\ref{bgd}), $P$-a.s.,
\begin{equation}
\label{bgdij} \beta_i(W^{i}_1)^2+\beta_j(W^{j}_1)^2 +\gamma_{ij}\int_0^1W^i_sdW^j_s+\gamma_{ji}\int_0^1W^j_sdW^i_s+\delta+\sum_{k\neq i,j}\beta^2_k \leq 0.
\end{equation}
Introducing the Levy area $L^{ij}=\int_0^1W^i_sdW^j_s-\int_0^1W^j_sdW^i_s$, we can write :
\[ \gamma_{ij}\int_0^1W^i_sdW^j_s+\gamma_{ji}\int_0^1W^j_sdW^i_s
= \frac 12(\gamma_{ij}+\gamma_{ji})W^i_1W^j_1+ \frac 12(\gamma_{ij}-\gamma_{ji})L^{ij}.
\]
If we substitute this in (\ref{bgdij}), it follows that, $P$-a.s.,
\[ \frac 12(\gamma_{ij}-\gamma_{ji})E[ L^{ij}| W^i_1=W^j_1=0]+\delta+\sum_{k\neq i,j}\beta^2_k \leq 0.\]
But, even after conditioning by $W^i_1=W^j_1=0$, the distribution of $L^{ij}$ is symmetric and
of unbounded support.
Consequently it holds that $\gamma_{ij}=\gamma_{ji}$. \\
Since  $(i,j)\in\{ 1,\ldots,d\}^2,\; i\neq j$ was chosen arbitrarily, (\ref{bgd}) becomes now, $P$-a.s.,
\[
 \sum_{i=1}^d\beta_i(W^{i}_1)^2+ \sum_{i< j}\gamma_{ij}W^i_1W^j_1+\delta \leq 0,
\]
or, equivalently,
\begin{equation}
\label{bgd3}
\frac 12 \langle W_1,AW_1\rangle+\delta\leq 0.\end{equation} 
Since the support of $W_1$ is $\R^ d$,
ii) and iii) follow.
\cq

\noindent {\bf Proof of the Theorem.}
We consider the following two additional assertions, where $u\in U$ is identified with the deterministic constant control process $u_t=u,t\geq 0$. Notice that $X^{x,u}$ is defined by (\ref{sds})  and $x^{x,t,u,v}$ by (\ref{de}).\\

f) For all $u\in U$, $x\in K$ and $t\geq 0$, $P[\xxu_t\in K]=1$.\\

g) For all $u\in U$, $x\in K$ and any admissible control $v\in \BR$, the function $\xxuv(t)$ takes its values in $K$.\\

\noindent The proof will be organized as follows :
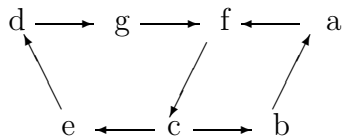
\begin{center}
\begin{picture}(100,50)
\put(20,0){e}
\put(60,0){c}
\put(100,0){b}
\put(0,40){d}
\put(40,40){g}
\put(80,40){f}
\put(120,40){a}
\put(55,2){\vector(-1,0){22}}
\put(70,2){\vector(1,0){22}}
\put(20,10){\vector(-1,2){14}}
\put(100,10){\vector(1,2){14}}
\put(10,42){\vector(1,0){22}}
\put(50,42){\vector(1,0){22}}
\put(110,42){\vector(-1,0){22}}
\put(75,35){\vector(-1,-2){14}}

\end{picture}
\end{center}
\begin{itemize}
\item d)$\Rightarrow$ g) is trivial.
\item g) $\Rightarrow$ f):
We fix some $u\in U$ and consider a scheme which converges in distribution
to the solution of the stochastic differential equation
(\ref{sds}) with constant, deterministic control $u_t=u$. For the construction we apply the following limit theorem \cite[Theorem 1, p. 698]{ls}:
For all $t\geq 0$, we set $\xi_t=W_{t+1}-W_t$. The process  $ {(\xi_t)}_{t \geq 0} $ is
strictly stationary and ergodic.
We let $ \eta_t^m = \sqrt{m} \xi_{mt} $, $ t \geq 0, m \geq 1 $, and put
$$
Y_t^m = \int_0^t \eta_s^m ds, t\geq 0.
$$
Notice that $ Y^m $ is a stochastic process with differentiable trajectories.
Furthermore,  the process $ Y^m$ converges $\omega$-wise, uniformly on compacts to $W$, as $ m \to \infty $. Consequently, it converges also in distribution
on pathspace. Theorem 1 from \cite{ls} tells now that the unique solution of
$$
dX_t^m = \tilde b(X_t^m,u)dt + \sigma(X_t^m,u)dY^m_t, \quad X_0^m = x,
$$
converges in distribution on pathspace to $ X^{x,u} $. The conditions
as stated in \cite{ls} on $ \sigma $ are slightly stronger than our
assumptions, namely $ C^2 $ is required. However, the proof in
\cite{ls} also holds for $\sigma$ satisfying our
$C^{1,1}$-assumptions. Certainly we cannot deduce by \cite[Theorem
1]{ls} a rate of convergence for $ X^m \to X^{x,u} $, but we also do
not need such a rate for our purposes.

We know that, by assumption, with probability 1, $ X_t^m \in K $ for all $ x \in K $, $ t \geq 0 $ and $ n \geq 1 $, whence we obtain the result: Indeed, if $d_K(x)$ denotes the distance from $x\in\R^n$ to $K$, we have
\[ E[d_K(X_t)]=\lim_{n\to \infty}E[d_K(X^m_t)]=0, \mbox{ for all }t\geq 0.\]
\item
f)$\Rightarrow$ c) :
Consider a constant control $u_t\equiv u\in U$ and suppose that, for all $x\in K$ and $t\geq 0$,
$P[X_t^{x,u }\in K]=1$. Let $\varphi\in C^2$ and
$x\in\argmax_K\varphi$. Up to change $\vf$ outside of some open
set including $x$, we can suppose that $\vf$, $\| D\vf\|$ and
$\|D^2\vf\|$ are bounded. We can apply the stochastic Taylor
expansion formula (\cite{lv} or \cite{fb}): for all $t\geq 0, P$-a.s.,
\[\begin{array}{rl}
 \varphi(\xxu_t)= &\varphi(x)+\sum_{i=1}^d\sigma^i_u\varphi(x)W^i_t+
\sum_{i=1}^d(\sigma^i_u)^2\varphi(x)\frac{(W^{i}_t)^2}2\\
\\
&+ \sum_{i\neq j}\sigma^j_u\sigma^i_u
\varphi(x)\int_0^tW^i_sdW^j_s +\tilde b_u
\varphi(x)t+R_t,
\end{array}\]
where $R_t$ satisfies $\frac{R_t}t\to 0$ in probability as $t\searrow 0$. We apply here the operator-notations
$\sigma^ i_u\vf(x)=\langle \sigma(x,u),D\vf(x)\rangle$ and 
$\tilde b_u
\varphi(x)=\langle \tilde b(x,u),D\vf(x)\rangle$.\\
Since $K$ is invariant for the constant control $u$ and since
$x\in\argmax_K\varphi$, we have $P$-a.s., for all $t\geq 0$, $
\varphi(\xxcu_t)\leq\varphi(x)$. Thus, $P$-a.s., for any fixed
$t\geq 0$,
\[\begin{array}{l}
\sum_{i=1}^d\sigma^i_u\varphi(x)W^i_t+
\sum_{i=1}^d(\sigma^i_u)^2\varphi(x)\frac{W^{i2}_t}2\\
\\
+ \sum_{i\neq j}\sigma^i_u\sigma^j_u \varphi(x)\int_0^tW^i_sdW^j_s
+\tilde b_u\varphi(x)t+R_t\leq 0.
\end{array}\]
Now we can apply Lemma \ref{lemmaabc} and get exactly the claim. 

\item c)$\Rightarrow$ b) becomes trivial as soon we write
\[ b_u\vf(x)+\frac 12\mbox{tr}(D^2\varphi(x)\sigma(x,u)\sigma^*(x,u))=\frac 12 A_{\vf,x}+\tilde b_u\vf(x).
\] 

\item b)$\Rightarrow$
a) : The proof is adapted from the equivalent result about
viability in \cite{bcq}. It is easy to see that,
if b) holds, then the map $f:x\mapsto 1-\i_K(x)$ is a
supersolution of
\[ \sup_{u\in U}\LR_{x,u}f(x)=0.\]
We consider now a constant $C\geq 1$
and an uniformly continuous application $g$ from $\R^n$ to $[0,1]$ that satisfies
$g(x)=0\;$ if and only if $\; x\in K$.
Since, for all $x\in \R^n$, $g(x)\leq Cf(x)$, $f$ is also a supersolution of
the following Hamilton-Jacobi-Bellman equation
\begin{equation}
\label{hjb}
\sup_{u\in U}\LR_{x,u}f(x)+g(x)-Cf(x)=0.
\end{equation}
But we know that the unique solution $V$ with polynomial growth of (\ref{hjb}) can be represented as
\[ V(x)=\sup_{u\in\UR}E[\int_0^\infty e^{-Cs}g(X^{x,u}_s)ds].\]
By the comparison theorem, we then  have
\[ V(x)\leq f(x), x\in\R^n.\]
For $x\in K$, this implies that, for all $u\in \UR$, for all
$t\geq 0$, $P[X^{x,u}_t\in K]=1$.

\item a) $\Rightarrow$ f) is trivial.

\item c)$\Rightarrow$ e) is trivial.

\item e)$\Rightarrow$ d) could be deduced from b)$\Rightarrow$ a)
if $v$ would take its values in a compact space and if $b$ and $\sigma
$ would be replaced by suitable functions. Let us clarify this
point: We fix $x \in K$,  $v \in \BR$ and $u \in \AR$. We wish to
prove that
$$x^{x,u,v}(t) \in K ,\mbox{ for all } t \geq 0.$$
 For any integer $n \geq 0$, we can define the control $t \mapsto v_n(t) :=
 \pi_n(v(t))$, where $\pi_n$ denotes the projection onto $\overline{B(0,n)}$.

By standard estimates, the sequence $x^{x,u,v_n} $ converges to
$x^{x,u,v}$ uniformly on every compact intervals  $[0,T]$.
Obviously  $x^{x,u,v_n} $ is solution to the following control
system
\begin{equation}\label{new}x'(t) = \widetilde{b} (x(t),u(t)) + \sigma
(x(t),u(t)) v(t) , u(t) \in U, v(t) \in \overline{B(0,n)}, x(0)=x,
\end{equation} with $(u,v)$ taking values in the compact set $U \times \overline{B(0,n)}$.
Hence we can apply the already proved relation b) $\Rightarrow$ a)
to the control system (\ref{new}) with $b (x,u)$ replaced by
$\widetilde{b}(x,u)  + \sigma (x,u) v$, $\sigma $ replaced by $0$
and the control $u$ replaced by $(u,v)$. In this case the relation
(\ref{5}) reduces to
$$\sup_{u\in \,U,\; v \in \, \overline{B(0,n)}) }\{ \LR'_{x,u}\varphi(x) +
\langle\sigma(x,u)v,D\varphi (x)\rangle \} \leq 0.$$ 
Consequently, we deduce from e)
that $x^{x,u,v_n}(t) \in K$, for all $t\geq 0$. By passing to the
limit with respect to  $n$, we obtain that $x^{x,u,v}(t) \in K $ forall $t \geq
0$. Our claim is proved.

\end{itemize}
\cq


\begin{thebibliography}{abcd99xyzu}

\bibitem{ada} AUBIN J.-P., DA PRATO G. (1998)
{\em The viability theorem for stochastic differential inclusions\/},
Stochastic Anal. Appl. 16, pp.1-15.

\bibitem{ad} AUBIN J.-P., DOSS H. (2003) 
{\em Characterization of stochastic viability of any nonsmooth set involving its generalized contingent curvature\/}, 
Stochastic Anal. Appl. 21, pp. 955-981.

\bibitem{bg} BARDI M. GOATIN P. (1999)
{\em Invariant sets for controlled degenerate diffusions: a viscosity solution approach\/}, 
in ``Stochastic Analysis, Control, Optimization and Applications: A Volume in Honor of W.H. Fleming'', W.M. McEneaney, C.G. Yin and Q. Zhang eds., Birkhaeuser, Boston, pp.191-208.

\bibitem{bj} BARDI M., JENSEN R. (2002)
{\em A geometric characterization of viable sets for controlled degenerate diffusions\/}, 
Set-Valued Anal. 10, pp. 129-141.

\bibitem{fb} BAUDOUIN F.(2004)
{\em An introduction to the Geometry of Stochastic Flows\/},
Imperial College Press.

\bibitem{bcq} BUCKDAHN R., CARDALIAGUET C., QUINCAMPOIX M. (2002)
{\em Representation Formula and mean Curvature Motion\/}, SIAM
J. Math. Anal. V. 33, N. 4, pp. 827-846.

\bibitem{bpqr} BUCKDAHN R., PENG S., QUINCAMPOIX M., RAINER C. (1998) 
{\em Existence of stochastic control under state constraints\/}, 
C.R.Acad.Sci. Paris S\'er. I Math. 327, pp. 17-22.

\bibitem{dpf} DA PRATO G., FRANKOVSKA H. (2001)
{\em Stochastic viability for compact sets  in terms of the distance function\/}, Dynamic Systems Appl. 10, pp. 177-184.

\bibitem{df} DA PRATO G., FRANKOWSKA H. (2004)
{\em Invariance of stochastic control systems with deterministic arguments\/},
Journal of Differential Equations 200, pp. 18-52.

\bibitem{d} DOSS H. (1977)
{\em Liens entre \'equations diff\'erentielles stochastiques et ordinaires\/},
 Ann. Inst. H. Pointcar\'e, Calcul Probab. Statist. 23, pp. 99-125.

\bibitem{gt} GAUTIER S., THIBAULT L. (1993)
{\em Viability for constrained stochastic differential equations\/},
Differential Integral Equations 6, pp. 1395-1414.

\bibitem{ls} LIPTSER R.Sh., SHIRYAYEV, A.N. (1989)
{\em Theory of Martingales\/}
Mathematics and its Applications
(Soviet Series), 49. Kluwer Academic Publishers Group.

\bibitem{lv} LYONS T., VICTOIR, N. (2004) {\em Curbature on Wiener Space\/},
Proceedings of the Royal Society London A 460, no. 2041, pp. 169-198.

\bibitem{m} MILIAN A. (1993) 
{\em A note on  stochastic invariance for It\^o equations\/}, Bull. Polish Acad. Sci. 41, pp. 139-150.

\bibitem{qr} QUINCAMPOIX M., RAINER C. (2005) {\em Stochastic control and compatible subsets of contraints\/}, Bull. Sci. math. 129, pp. 39-55.

\bibitem{sus} SUSSMANN H.J. (1978)
{\em On the gap between deterministic and stochastic ordinary differential equations\/},
Ann. Probability 6, pp. 19-41.

\end{thebibliography}
\end{document}